\documentclass[reqno]{amsart}

\usepackage{amsfonts}
\usepackage{amssymb}
\usepackage{tikz}
\usepackage{hyperref}

\copyrightinfo{2010}{American Mathematical Society}

\newtheorem{theorem}{Theorem}[section]
\newtheorem{lemma}[theorem]{Lemma}

\theoremstyle{definition}
\newtheorem{definition}[theorem]{Definition}
\newtheorem{example}[theorem]{Example}

\theoremstyle{remark}

\numberwithin{equation}{section}

\newcommand{\CN}{\mathbb{C}}
\newcommand{\cD}{\mathcal{D}}

\newcommand{\R}{\mathbb{R}}

\newcommand{\cC}{\mathcal{C}}
\newcommand{\cA}{\mathcal{A}}

\begin{document}

\title[A Class of Toeplitz Operators in Several Complex Variables]
      {A Class of Toeplitz Operators in Several Complex Variables}

\author{D. Fedchenko}

\address[Dmitry Fedchenko]
        {Institute of Mathematics and Computer Science,
         Siberian Federal University,
         Svobodny Prospekt 79,
         660041 Krasnoyarsk,
         Russia}

\email{fdp@bk.ru}





\date{November 10, 2016}




\begin{abstract}
In order to study the Toeplitz algebras related to a Dirac operators in a neighborhood of a closed bounded domain $\cD$ with smooth boundary in $\CN^n$ we introduce a singular Cauchy type integral. We compute its principal symbol, thus initiating the index theory.
\end{abstract}

\maketitle

\tableofcontents

\section*{Introduction}
\label{s.Introduction}
There are a number of ways in which the theory of Toeplitz operators can be generalised
to $n$ dimensions, see e.g.
   \cite{Venu72},
   \cite{Doug73}
and the references given there.
The monograph \cite{Upme96} presents much more advanced theory of Toeplitz operators in
several complex variables.

The paper \cite{Guil84} describes precisely how Toeplitz operators of ``Bergman type''
are related to Toeplitz operators of ``Szeg\"{o} type''.
A remarkable connection between the theory of Toeplitz operators \'{a} la \cite{Venu72}
and the standard theory of pseudodifferential operators emerged from the work
   \cite{BoutGuil81}.
This connection in its broad outlines is elucidated in \cite{Guil84}, too.

This work focuses on a new class of Toeplitz operators which is more closely related to
elliptic theory.
The new Toeplitz operators admit very transparent description which motivates strikingly
their study.
To this end, let $\cA$ be a $(k \times k)\,$-matrix of first order scalar partial differential operators
in a neighborhood of the closed bounded domain $\cD$ with smooth boundary $S$ in $\CN^n$.
Assume that the leading symbol of $\cA$ has rank $k$
away from the zero section of the cotangent bundle of $\cD$.
Then, given any solution $u$ of the homogeneous equation $\cA u = 0$ in the interior of $\cD$
which has finite order of growth at the boundary, the Cauchy data $t (u)$ of $u$ with respect
to $\cA$ possess weak limit values at the boundary.
If $\cA$ satisfies the so-called uniqueness condition of the local Cauchy problem in a
neighborhood of $\cD$, then the solution $u$ is uniquely defined by its Cauchy data at
$S$.
Let
   $t (u) = B u$
be a representation of the Cauchy data of $u$.
The space of all Cauchy data $B u$ of $u$ at the boundary is effectively
described by the condition of orthogonality to solutions of the formal adjoint equation
$\cA^* g = 0$ near $\cD$ by means of a Green formula, see \cite[\S~10.3.4]{Tar95}.
In this way we distinguish Hilbert space of vector-valued functions on $S$ which represent solutions to $Au = 0$ in the interior of $\cD$.
In particular, one introduces Hardy spaces $H$ as subspaces of
   $L^2 (S, \CN^k)$
consisting of the Cauchy data of solutions to $\cA u = 0$ in the interior of $\cD$ with
appropriate behaviour at the boundary.
Pick such a Hilbert space $H$.
By the above, $H$ is a closed subspace of $L^2 (S,\CN^k)$ and we write
$\Pi$ for the orthogonal projection of $L^2 (S,\CN^k)$ onto $H$. If $\cA$ is the Cauchy--Riemann operator then $\Pi$ just amounts to the Szeg\"{o} projection.

Given a $(k \times k)\,$-matrix $M (z)$ of smooth function on $S$, the operator
$T_M$ on $H$ given by $u \mapsto \Pi (Mu)$ is said to be a Toeplitz operator with
multiplier $M$.
If $M$ is a scalar multiple of the unit matrix, then the theory of Toeplitz operators $T_M$
is much about the same as the classical theory of Toeplitz operators. Otherwise the theory is much more complicated.

Finally, we briefly discuss Toeplitz operators related to the algebra of octonions $\mathbb{O}$ introduced by John T. Graves in 1843. The octonions were discovered independently by Cayley and are sometimes referred to as Cayley numbers or the Cayley algebra.

\section{Dirac operators}
\label{s.Dirac operators}
Let $\CN^n$ be the standard $n$-dimensional complex space obtained from the underlying real space $\R^{2n}$ of variables $x=(x_1, \dots, x_{2n})$ by introducing the complex structure
$$
z_j = x_j + \imath x_{n+j}
$$
for $j=1, \dots, n$. As usual, we define complex derivatives by
\begin{align*}
\frac{\partial}{\partial z_j} &= \frac{1}{2} \left( \frac{\partial}{\partial x_j} - \imath \frac{\partial}{\partial x_{n+j}} \right),\\
\frac{\partial}{\partial \overline z_j} &= \frac{1}{2} \left( \frac{\partial}{\partial x_j} + \imath \frac{\partial}{\partial x_{n+j}} \right).
\end{align*}

Let us consider a matrix-valued constant coefficient operator
\begin{equation}\label{DiracOp}
\cA \left( \frac{\partial}{\partial z}, \frac{\partial}{\partial \overline z} \right) = \sum_{j=1}^n \alpha_j \frac{\partial}{\partial z_j} + \beta_j \frac{\partial}{\partial \overline z_j}
\end{equation}
such that
\begin{equation}\label{Laplace_Factor}
\cA^* \cA = - \frac{1}{4} \Delta E
\end{equation}
where $\cA^*$ is the formal adjoint of $\cA$ and $E$ the identity matrix. It is easily seen that the identity (\ref{Laplace_Factor}) reduces to a system of identities for the coefficients, namely
\begin{align}\label{Dirac_Coef}
\alpha_j^* \alpha_k + \beta^*_k \beta_j &= \delta_{j,k} E, \notag\\
\alpha^*_j \beta_k + \alpha^*_k \beta_j &= 0
\end{align}
for all $j,k = 1, \dots, n$.

It is well known that there is a solution of (\ref{Dirac_Coef}) amongst matrices of type $(2^{n-1} \times 2^{n-1})$, cf. for instance Chapter 3 in \cite{Me93}.

First-order differential operators with constant coefficients factorizing the Laplacian in the sense of (\ref{Laplace_Factor}) are called Dirac operators.

\begin{example}
The Cauchy--Riemann operator $\cA = \partial / \partial \bar{z}$ is a Dirac operator in the complex plane.
\end{example}

\begin{example}
The operator
$$
\cA
 = \left( \begin{array}{rr}
              \partial / \partial \bar{z}_1
          & - \partial / \partial \bar{z}_2
          \\
              \partial / \partial z_2
          &   \partial / \partial z_1
          \end{array}
   \right)
$$
is a Dirac operator in $\mathbb{C}^2$.
\end{example}

\section{Fundamental solution}
To construct an explicit fundamental solution $\Phi$ for a Dirac operator $\cA$ we make use of relation (\ref{Laplace_Factor}). Namely, denote by $e$ the standard fundamental solution of convolution type for the Laplace operator in $\R^{2n}$. In the coordinates of $\CN^n$ it reads
$$
e(z) = \frac{(n-1)!}{2 \pi^n} \frac{1}{2-2n} \frac{1}{|z|^{2n-2}},
$$
if $n>1$, and $e(z) = (1/2\pi) \ln |z|$, if $n=1$.

Set
\begin{align*}
\Phi(z - \zeta) &= 4 \cA^*\left( \frac{\partial}{\partial \zeta}, \frac{\partial}{\partial \overline \zeta} \right) e(z-\zeta)
\\
&= \frac{(n-1)!}{\pi^n} \frac{\cA^*(\overline z - \overline \zeta, z - \zeta)}{|z-\zeta|^{2n}}
\end{align*}
for $z \neq \zeta$, where by $\cA^*(\overline z - \overline \zeta, z - \zeta)$ is meant the adjoint of the matrix $\cA(\overline z - \overline \zeta, z - \zeta)$.

\begin{lemma}\label{Fund_Sol}
As defined above, $\Phi(z-\zeta)$ is a fundamental solution of the Dirac operator $\cA$, i.e. $\Phi \circ \cA = E$ and $\cA \circ \Phi = E$ on $C_{\rm comp}^\infty(\CN^n, \CN^k)$, where $(k \times k)$ is the type of $E$.
\end{lemma}

{\it Proof.} The first relation $\Phi \circ \cA = E$ is fulfilled by the very construction of $\Phi$. Since $\cA$ is a square matrix, if follows from (\ref{Laplace_Factor}) that $\cA \cA^* = -(1/4)\Delta E$ whence $\Phi^* \circ \cA^* = E$ on $C_{\rm comp}^\infty(\CN^n,\CN^k)$. The latter equality is equivalent to $\cA \circ \Phi = E$, as desired.
\hfill $\square$

\section{Green formula}
Let $\cD$ be a bounded domain with smooth boundary $S = \partial \cD$ in $\CN^n$. Write $\nu(y) = (\nu_1(y), \dots, \nu_{2n}(y))$ for the unit outward normal vector of $S$ at a point $y \in S$.

If $\rho(x)$ is a defining function of $S$ then
$$
\nu (y) = \frac{\nabla \rho (y)}{|\nabla \rho (y)|}
$$
for $y \in S$, where $\nabla \rho (y)$ stands for the real gradient of $\rho$ at $y$. The complex vector $\nu_c = (\nu_{c,1}, \dots, \nu_{c,n})$ with coordinates $\nu_{c,j} = \nu_j + \imath \nu_{n+j}$ is called the complex normal of the hypersurface $S$. In the coordinates of $\CN^n$ we obviously have
$$
\nu_{c,j} = \frac{\partial \rho / \partial \overline \zeta_j}{|\nabla_{\overline \zeta} \rho(\zeta)|}
$$
for $j=1, \dots, n$.

Denote by $d \zeta$ the wedge product $d \zeta_1 \wedge \dots \wedge d \zeta_n$, and by $d \overline \zeta [j]$ the wedge product of all $d \overline \zeta_1, \dots, d \overline \zeta_n$ but $d \zeta_j$.

\begin{lemma}\label{pullback}
For each $j=1,\dots,n$, the pullback of the differential form $d \zeta \wedge d \overline \zeta[j]$ under the embedding $S \hookrightarrow \CN^n$ is equal to $(-1)^{j-1}(2 \imath)^{n-1} \imath \nu_{c,j} ds$, where $ds$ is the area form on the hypersurface $S$ induced by the Hermitian metric of $\CN^n$.
\end{lemma}

{\it Proof.} An easy computation shows that the pullback of the differential form $dy[j]$ under the embedding $S \hookrightarrow \CN^n$ is equal to $(-1)^{j-1} \nu_j ds$, for every $j=1, \dots, 2n$. From this the lemma follows immediately.
\hfill $\square$

Denote by $\sigma(\zeta)$ the principal homogeneous symbol $\sigma^1(\cA)(\zeta, \xi)$ of the operator $\cA$ evaluated at the point $(\zeta, \nu_c(\zeta) / \imath)$ of the cotangent bundle of a neighborhood of $\cD$, where $\zeta \in S$. For example, the principal homogeneous symbol of the Cauchy--Riemann operator is equal $-(1/2)(\xi_1 + \imath \xi_2)$.

We are now in a position to specify the restriction of the Green operator $G_\cA(g,u)$ of $\cA$ to the boundary. By a Green operator of $\cA$ is meant a bilinear operator $G_\cA$ from $C^\infty(\CN^n, (\CN^k)^*) \times C^\infty(\CN^n, \CN^k)$ to differential forms of degree $2n-1$ on $\CN^n$, such that $dG_\cA(g^*,u) = ((\cA u,g)_y - (u, \cA^*g)_y)dy$ holds pointwise in $\CN^n$.

By \cite[$\S$ 2.4.2]{Tar90}, there is a unique Green operator for $\cA$, and its pullback under the embedding $S \hookrightarrow \CN^n$ is

\begin{align*}
G_\cA(g,u) &= \frac{1}{\imath} g \cA \left( \frac{\imath}{2} (\nu_j - \imath \nu_{n+j}), \frac{\imath}{2}(\nu_j + \imath \nu_{n+j}) \right) u ds\\
         &= g \sigma(\zeta) u ds
\end{align*}
whence
\begin{equation}\label{GreenOp}
G_\cA (\Phi(z-\zeta), E) = \frac{(n-1)!}{\pi^n} \frac{\cA^*(\overline z - \overline \zeta, z-\zeta)}{|z-\zeta|^{2n}} \sigma(\zeta) ds.
\end{equation}

\begin{lemma}\label{GreenFormula}
Every vector-valued function $u \in C^1(\overline \cD, \CN^k)$ has the integral representation
$$
\chi_\cD u = - \int_S G_\cA (\Phi(z-.), u) + \int_\cD \Phi(z-.) \cA u \, dv,
$$
where $dv$ is the Lebesgue measure in $\R^{2n}$ and $\chi_\cD$ the characteristic function of $\cD$.
\end{lemma}

{\it Proof.} This is a very special case of a general Green formula related to an elliptic system of differential equations, see for instance \cite[$\S$ 2.5.4]{Tar90} and elsewhere.
\hfill $\square$

Needless to say that this formula extends to the case of Sobolev class functions $u \in H^1(\cD, \CN^k)$ as well as to more general functions on $\cD$.

\section{Toeplitz operators}\label{secTO}
\begin{theorem}\label{OrtProp}
Let $u_0 \in L^2(S, \CN^k)$. In order that there be a solution $u$ to $\cA u = 0$ in the interior of $\cD$, which has finite order of growth at $S$ and coincides with $u_0$ on $S$, it is necessary and sufficient that
\begin{equation}\label{IntOrtProp}
\int_S g \sigma(\zeta) u_0 \,ds = 0
\end{equation}
for all solutions of the formal adjoint equation $\cA^* g = 0$ near $\cD$.
\end{theorem}

{\it Proof.} See Theorem 10.3.14 in \cite{Tar95}.
\hfill $\square$

We denote by $H$ the (closed) subspace of $L^2(S, \CN^k)$ that consists of all functions $u$ satisfying the orthogonality conditions (\ref{IntOrtProp}). The elements of $H$ can be actually specified as solutions to $\cA u =0$ of Hardy class $H^2$ in the interior of $\cD$, see \cite[$\S$ 11.2.2]{Tar95}. The orthogonal projection $\Pi$ of $L^2(S, \CN^k)$ onto $H$ is therefore an analogue of Szeg\"o projection.

\begin{definition}\label{defTO}
Let $M(z)$ be a $(k \times k)$-matrix whose entries are bounded functions on $S$. By a Toeplitz operator $T_M$ with multiplier $M$ is meant the operator $u \mapsto \Pi(Mu)$ in $H$.
\end{definition}

More generally, if $\Psi$ is a $(k\times k)$-matrix of pseudodifferential operators of order $0$ on $S$, then $\Psi$ maps $L^2(S, \CN^k)$ continuously into $L^2(S, \CN^k)$. Therefore, the composition $T_\Psi = \Pi \Psi$ is a continuous self-mapping of $H$ which we call a generalised Toeplitz operator.

If the projection $\Pi$ is a classical pseudodifferential operator on $S$, then from the equality $\Pi^2 = \Pi$ it follows readily that the order of $\Pi$ just amounts to zero. Hence, the generalised Toeplitz operators on $S$ form a subalgebra of the $C^*$-algebra of all zero order classical pseudodifferential operators on $C^\infty(S,\CN^k)$.

\section{The generalised Cauchy type integral}

To clarify the nature of the generalised Szeg\"o projection $\Pi$ we introduce the singular Cauchy type integral
\begin{equation}\label{SCI}
\cC u(z) = - \textrm{p.v.} \int_S G_\cA (\Phi(z-.), u)
\end{equation}
for $z \in S$, where $u \in L^2(S, \CN^k)$. The principal value of the integral on the right-hand side exists for almost all $z \in S$ and it induces a bounded linear operator in $L^2(S, \CN^k)$.

\begin{lemma}\label{PBformula}
The operators $(1/2)I \pm \cC$ are projections on the space $L^2(S, \CN^k)$.
\end{lemma}

{\it Proof.} This follows from the equality $\cC^2 = (1/4)I$ by a trivial verification, cf. for instance \cite{Tar06}.
\hfill $\square$

The generalised Cauchy type integral (\ref{SCI}) is a classical pseudodifferential operator of order 0 in $C^\infty(S, \CN^k)$. We finish this section by evaluating its principal homogeneous symbol. To this end, we identify the cotangent space $T^*_z S$ of $S$ at a point $z \in S$ with all linear forms on $T^*_z \R^{2n}$ which vanish on the one-dimensional subspace of $T^*_z \R^{2n}$ spanned by $\nu(z)$. Since $T^*_z \R^{2n} \cong \R^{2n}$, one can actually specify $T^*_z S$ as the hyperplane through the origin in $\R^{2n}$ which is orthogonal to the vector $\nu(z)$.

\begin{lemma}\label{SymbEval}
For each $z \in S$ and $\xi \in T^*_z S$, the symbol of order 0 of the operator $\cC$ is equal to
$$
\sigma^0(\cC)(z,\xi) = \sigma^1(\cA^*) \left(z, \frac{\xi}{|\xi|}\right) \, \sigma(z).
$$
\end{lemma}

{\it Proof.} By assumption, the Laplacian $\cA^* \cA$ is a second order elliptic differential operator on $C^\infty(U, \CN^k)$. It has a parametrix $P$ which is a $(k \times k)$-matrix of scalar pseudodifferential operators of order $-2$ on $U$. The operator $\Phi$ differs from $P \cA^*$ by a smoothing operator, and so it has the principal homogeneous symbol $(\sigma^2(\cA^* \cA))^{-1} \sigma^1(\cA^*)$ which is a left inverse for $\sigma^1(\cA)$.

Formula (\ref{SCI}) just amounts to saying that
$$
\cC u(z) = - \Phi(\sigma(\zeta) u(\zeta) \ell_S),
$$
where $\ell_S$ is the surface layer on $S$. We thus see that the pseudodifferential operator $\cC$ on $S$ is the restriction to $S$ of the pseudodifferential operator
$$
\Psi = - \Phi \circ \sigma(z)
$$
defined in a neighborhood of $S$. This latter is of order $-1$ and its principal symbol is easily evaluated, namely
\begin{align*}
\sigma^{-1}(\Psi)(z,\xi) &= - \sigma^{-1}(\Phi)(z, \xi) \, \sigma(z)\\
&= -(\sigma^2(\cA^* \cA)(z,\xi))^{-1} \, \sigma(z)
\end{align*}
for $z$ in a neighborhood of $S$ and $\xi \in \CN^n \setminus \{0\}$. A familiar argument now shows that the principal symbol of $\cC$ is given by the formula
\begin{align*}
\sigma^0(\cC)(z,\xi) &= \frac{1}{2 \pi} {\rm p.v.} \int_{-\infty}^\infty \sigma^{-1}(\Psi)(z, t \nu(z)+\xi)\, dt\\
&= -\frac{1}{2 \pi} {\rm p.v.} \int_{-\infty}^\infty (\sigma^2(\cA^* \cA)(z,t \nu(z) + \xi))^{-1} \, \sigma^1(\cA^*)(z,t \nu(z) + \xi)\, dt\\
&\times \sigma(z)
\end{align*}
for all $z \in S$ and $\xi \in \R^{2n} \setminus \{0\}$ orthogonal to $\nu(z)$. Note that the integral on the right-hand side diverges, however, its Cauchy principal value exists, which is due to the condition $\langle \nu(z), \xi \rangle =0$.

Since
$$
\sigma^2(\cA^* \cA)(z, \xi) = \frac{1}{4}|\xi|^2 E,
$$
we shall have established the desired equality if we prove that
$$
\frac{1}{4 \pi} \int_{-\infty}^\infty \frac{1}{|t \nu(z) + \xi|^2} \, dt = \frac{1}{4} \frac{1}{|\xi|}
$$
for all $z \in S$ and $\xi \in \R^{2n} \setminus \{0\}$ orthogonal to the vector $\nu(z)$. A trivial verification shows that
$$
|t \nu(z) + \xi|^2 = t^2 + |\xi|^2
$$
whence
\begin{align*}
\frac{1}{4 \pi} \int_{-\infty}^\infty \frac{1}{|t \nu(z)+\xi|^2} \, dt &= \frac{1}{4 \pi} \int_{-\infty}^\infty \frac{1}{t^2 + |\xi|^2} \, dt
\\
&= \frac{1}{4} \frac{1}{|\xi|},
\end{align*}
as desired.
\hfill $\square$

\section{Index of Toeplitz operators}
We now return to generalised Toeplitz operators $T_\Psi = \Pi \Psi$ in $H$ introduced in Section \ref{secTO}, see Definition \ref{defTO} and below. Here, $\Pi$ is a projection of $L^2(S, \CN^k)$ onto $H$. We restrict ourselves to the case where $\Pi$ is a classical pseudodifferential operator of order zero in $C^\infty(S, \CN^k)$. We are interested in characterizing those Toeplitz operators in $H$ which possess the Fredholm property. To this end we extend $T_\Psi$ from $H$ to all of $L^2(S, \CN^k)$ in a special manner and use Fredholm criteria for operator algebras with symbols.

\begin{lemma}
A Toeplitz operator $T_\Psi$ in $H$ is Fredholm if and only if so is the operator
$$
E_\Psi := T_\Psi \Pi + (I - \Pi)
$$
in $L^2(S, \CN^k)$.
\end{lemma}

{\it Proof.} We see that $T_\Psi$ is the restriction of the pseudodifferential operator
$$
E_\Psi = T_\Psi \Pi \oplus (I-\Pi): H \oplus H^\bot \to H \oplus H^\bot
$$
on $H$.

If $E_\Psi$ is Fredholm then ${\rm ind \,} E_\Psi$ is finite. But ${\rm ind \,} E_\Psi = {\rm ind \,} T_\Psi + {\rm ind \,} (I-\Pi)$, where $T_\Psi : H \to H$ and $(I-\Pi): H^\bot \to H^\bot$. It is clear that ${\rm ind \,} (I-\Pi)=0$, whence ${\rm ind \,} T_\Psi$ is also finite. And moreover, ${\rm ind \,} E_\Psi = {\rm ind \,} T_\Psi$.
\hfill $\square$

The operator $E_\Psi$ on $L^2(S, \CN^k)$ is pseudodifferential of order zero. Its principal homogeneous symbol just amounts to $\sigma^0(E_\Psi) = \sigma^0(\Pi) \sigma^0(\Psi) \sigma^0(\Pi) + \sigma^0(I - \Pi)$ away from the zero section of $T^* S$. Given any $s \in \R$, the operator $E_\Psi$ in $H^s(S, \CN^k)$ is known to be Fredholm if and only if it is elliptic, i.e. $\sigma^0(E_\Psi)(z,\xi)$ is invertible for all $(z, \xi) \in T^* S$ with $\xi \neq 0$. Moreover, the index of this operator is actually independent of the particular choice of $s$ and it can be evaluated by the familiar Atiyah--Singer formula \cite{AS68}.



\section{Concluding remarks}
In the sequel we restrict our attention to a special Clifford algebra corresponding to the case $n=4$. Clifford algebras have important applications in variety of fields including string theory, special relativity and quantum logic. The algebra in question is called the algebra of octonions and denoted by $\mathbb{O}$. To wit,
$$
\cA
 = \left( \begin{array}{rr}
              q_1
          & - q_2
          \\
              \bar{q}_2
          &   \bar{q}_1
          \end{array}
   \right),
$$
where $q_i = \left( \frac{\partial}{\partial \bar{z}_i} + \frac{\partial}{\partial \bar{z}_{2+i}} \jmath \right)$, $\jmath$ is the fundamental quaternion unit and the bar denotes the conjugate of a quaternion.

Let us rewrite the action of the operator $\cA$ on a quaternion-valued function in complex coordinates, that is
\begin{align*}
\cA u
 &= \left( \begin{array}{rr}
              \bar{\partial}_1 + \bar{\partial}_3 \jmath
          & - \bar{\partial}_2 - \bar{\partial}_4 \jmath
          \\
              \partial_2 - \bar{\partial}_4 \jmath
          &   \partial_1 - \bar{\partial}_3 \jmath
          \end{array}
   \right)\,
\left( \begin{array}{r}
              u_1 + u_3 \jmath
          \\
              u_2 + u_4 \jmath
          \end{array}
   \right)
\\
&= \left( \begin{array}{rrrr}
              \bar{\partial}_1 \phantom{\mathfrak{c}}
          & - \bar{\partial}_2 \phantom{\mathfrak{c}}
          & - \bar{\partial}_3 \mathfrak{c}
          &   \bar{\partial}_4 \mathfrak{c}
          \\
              \partial_2 \phantom{\mathfrak{c}}
          &   \partial_1 \phantom{\mathfrak{c}}
          &   \bar{\partial}_4 \mathfrak{c}
          &   \bar{\partial}_3 \mathfrak{c}
          \\
              \bar{\partial}_3 \mathfrak{c}
          & - \bar{\partial}_4 \mathfrak{c}
          &   \bar{\partial}_1 \phantom{\mathfrak{c}}
          & - \bar{\partial}_2 \phantom{\mathfrak{c}}
          \\
            - \bar{\partial}_4 \mathfrak{c}
          & - \bar{\partial}_3 \mathfrak{c}
          &   \partial_2 \phantom{\mathfrak{c}}
          &   \partial_1 \phantom{\mathfrak{c}}
          \end{array}
   \right)\,
\left( \begin{array}{r}
              u_1
          \\
              u_2
          \\
              u_3
          \\
              u_4
          \end{array}
   \right),
\end{align*}
$\mathfrak{c}$ being the complex conjugation and $\bar{\partial}_j = \partial / \partial \bar{z}_j$.

The octonions have the dimension eight. Because of nonassociativity they can not be represented as quaternion $(2 \times 2)$-matrices with usual multiplication. The product of two octonions $\mathfrak{O}=(a,b)$ and $\mathfrak{P}=(c,d)$ is defined via the Cayley--Dickson construction by
$$
(a,b)(c,d) = (ac-\bar{d} b, da + b\bar{c}).
$$
It corresponds to the special multiplication of the matrices
$$
\mathfrak{O} \mathfrak{P} =
\left( \begin{array}{rr}
              a
          & - b
          \\
              \bar{b}
          &   \bar{a}
          \end{array}
   \right)
\,
\left( \begin{array}{rr}
              c
          & - d
          \\
              \bar{d}
          &   \bar{c}
          \end{array}
   \right)
=
\left( \begin{array}{rr}
              ac - \bar{d}b
          & -da - b \bar{c}
          \\
              c \bar{b} + \bar{a} \bar{d}
          &   -\bar{b}d + \bar{c} \bar{a}
          \end{array}
   \right).
$$

It remains to show that the multiplication $\mathfrak{O} \mathfrak{P}$ is alternative, i.e. $\mathfrak{O}(\mathfrak{O}\mathfrak{P})=(\mathfrak{O}\mathfrak{O})\mathfrak{P}$ and $(\mathfrak{P}\mathfrak{O})\mathfrak{O} = \mathfrak{P}(\mathfrak{O}\mathfrak{O})$ for all matrices $\mathfrak{O}$, $\mathfrak{P}$.

\begin{lemma}\label{MultAlt}
The matrices $\mathfrak{O}$ constitute an alternative algebra.
\end{lemma}

{\it Proof.}
\begin{align*}
\mathfrak{O} (\mathfrak{O} \mathfrak{P}) &=
\left( \begin{array}{rr}
              a & -b \\
              \bar{b} & \bar{a}
          \end{array}
   \right)
\,
\left( \begin{array}{rr}
              ac - \bar{d}b & -da - b\bar{c} \\
              c\bar{b} + \bar{a}\bar{d} & -\bar{b} d + \bar{c} \bar{a}
          \end{array}
   \right)
\\
&=
\left( \begin{array}{ll}
              aac - a \bar{d} b - c \bar{b} b - \bar{a}\bar{d} b & -daa-b\bar{c}a+b\bar{b}d-b\bar{c}\bar{a} \\
              ac\bar{b} - \bar{d} b \bar{b} + \bar{a} c \bar{b} + \bar{a}\bar{a}\bar{d} & -\bar{b}da-\bar{b}b\bar{c}-\bar{b}d\bar{a}+\bar{c}\bar{a}\bar{a}
          \end{array}
   \right)
\end{align*}
and
\begin{align*}
(\mathfrak{O} \mathfrak{O}) \mathfrak{P} &=
\left( \begin{array}{rr}
              aa-\bar{b}b & -ba-b\bar{a} \\
              a\bar{b}+\bar{a}\bar{b} & -\bar{b}b+\bar{a}\bar{a}
          \end{array}
   \right)
\,
\left( \begin{array}{rr}
              c & -d \\
              \bar{d} & \bar{c}
          \end{array}
   \right)
\\
&=
\left( \begin{array}{ll}
              aac - \bar{b}bc - \bar{d}ba - \bar{d}b\bar{a} & -daa+d\bar{b}b-ba\bar{c}-b\bar{a}\bar{c} \\
              ca\bar{b} + c\bar{a}\bar{b} - \bar{b}b\bar{d} + \bar{a}\bar{a}\bar{d} & -a\bar{b}d-\bar{a}\bar{b}d-\bar{c}\bar{b}b+\bar{c}\bar{a}\bar{a}
          \end{array}
   \right),
\end{align*}
where
\begin{align*}
a\bar{d}b+\bar{a}\bar{d}b-\bar{d}ba-\bar{d}b\bar{a} &= [a,\bar{d}b]+[\bar{a},\bar{d}b]
\\
&=[a+\bar{a},\bar{d}b]
\\
&=0
\end{align*}
and
\begin{align*}
ac\bar{b}+\bar{a}c\bar{b}-ca\bar{b}-c\bar{a}\bar{b} &= [\bar{a},c]\bar{b}+[a,c]\bar{b}
\\
&=[\bar{a}+a,c]\bar{b}
\\
&=0.
\end{align*}
The equality $(\mathfrak{P}\mathfrak{O})\mathfrak{O} = \mathfrak{P}(\mathfrak{O}\mathfrak{O})$ is verified similarly.
\hfill $\square$

To specify the class of multipliers we describe those $(2 \times 2)$-matrices $M(z)$ which commute with $\sigma^1(\cA)(z,\zeta)$.

\begin{lemma}\label{lemma_OM}
A $(2\times 2)$-matrix $M$ commutes with $\sigma^1(\cA)(z, \zeta)$ if and only if $M$ is of the form
\begin{equation}\label{Equ_OM}
M = \left( \begin{array}{ll}
              X
          & - YC
          \\
              YC
          &   \phantom{-} X
          \end{array}
   \right),
\end{equation}
C being quaternion conjugation and $X$, $Y$ are arbitrary $(2 \times 2)$-matrices with complex entries.
\end{lemma}

{\it Proof.} Here the multiplication $\mathfrak{O} M$ is considered in the usual way. Straightforward calculation shows that
$$
[\mathfrak{O},M] =
\left( \begin{array}{rr}
              [a,X]-[b,Y]C & -[b,X]-[a,Y]C \\
              \phantom{-} [\bar{b},X]+[\bar{a},Y]C & [\bar{a},X]-[\bar{b},Y]C
          \end{array}
   \right).
$$
An easy computation shows that $[a,X]=0$ if and only if $X$ is of the form
\begin{equation}\label{complex_mult}
X
 = \left( \begin{array}{ll}
              w_1
          & - w_2 \mathfrak{c}
          \\
              w_2 \mathfrak{c}
          &   \phantom{-} w_1
          \end{array}
   \right),
\end{equation}
where $w_1$, $w_2$ are arbitrary complex numbers.
\hfill $\square$

Needless to say that $X$ is a $(2 \times 2)$-matrix of nonlinear operators. From the equality

$$
X^1 X^2 =
\left( \begin{array}{rr}
              w^1_1 w^2_1 - w^1_2 \overline{w}^2_2 \phantom{) c}
          & - (w_2^1 \overline{w}_1^2 + w^1_1 w^2_2) \mathfrak{c}
          \\
              (w_2^1 \overline{w}_1^2 + w^1_1 w^2_2) \mathfrak{c}
          &   \phantom{- (} w^1_1 w^2_1 - w^1_2 \overline{w}^2_2 \phantom{) c}
          \end{array}
   \right)
$$
it follows readily that the matrices $X$ of the form (\ref{complex_mult}) constitute an unital algebra. Since $\det X = w_1^2 + |w_2|^2$, a matrix $X$ is invertible if and only if $X \neq 0$.

Since the kernel and cokernel of a nonlinear operator $X$ fail to be vector spaces, the index of $X$ is no longer
defined.
By the index of a differentiable nonlinear operator is usually meant the index of its Fr\'echet derivative.
However, the operators in question are not Fr\'{e}chet differentiable in the sense of complex vector spaces.

\begin{lemma}
The corresponding operator (\ref{complex_mult}) is not differentiable unless $w_2 \neq 0$.
\end{lemma}

{\it Proof.}
This follows by immediate computation.
\hfill $\square$

We conclude that there is no suitable Fredholm problems over the field of complex numbers for Toeplitz operators with multipliers of special kind.
From now on all consideration are over the field of real numbers.

Using Lemma \ref{PBformula} we establish a very useful formula for the projection $\Pi$, namely, $\Pi = (1/2)I + \cC$. Thus,
\begin{align*}
T_{M^1} = (1/2)M^1 + \cC M^1,\\
T_{M^2} = (1/2)M^2 + \cC M^2,
\end{align*}
and so
\begin{equation}\label{Semicommutator}
T_{M^1} T_{M^2} = T_{M^1 M^2} - [\cC, M^2]((1/2)M^1 - \cC M^1).
\end{equation}
\begin{lemma}\label{LemmaSC}
Assume that $M^1$ and $M^2$ are $(2 \times 2)$-matrices of smooth functions on $S$. Then $T_{M^2 M^1} = T_{M^2} T_{M^1}$ modulo compact operators on $H$.
\end{lemma}

{\it Proof.} If a $(2 \times 2)$-matrix $M$ commutes with $\sigma^1(\cA)(z,\xi)$, it is of the form (\ref{Equ_OM}), and so the adjoint matrix $M^*$ has the same form. Hence it follows that $M^*$ commutes with $\sigma^1(\cA)(z,\xi)$, and so $M$ commutes with the adjoin $\sigma^1(\cA^*)(z, \xi)$, too. Now an elementary analysis shows readily that
\begin{align*}
\cC(Mu)(z) &= - {\rm p.v.} \int_S M(\zeta) \Phi(z-\zeta) \sigma(\zeta) u(\zeta) ds\\
           &= M(\cC u)(z) + \int_S G_{\cA}(\Phi(z-\zeta), (M(z)-M(\zeta))u)_\zeta
\end{align*}
holds almost everywhere on the boundary for all $u \in L^2(S,\CN^4)$. In particular, if $u \in H$, then
$$
\cC(Mu)(z) = (1/2)(Mu)(z) + \int_S G_{\cA}(\Phi(z-\zeta), (M(z)-M(\zeta))u)_\zeta
$$
for almost all $z \in S$. From these two equalities we conclude that the remainder $[\cC, M^2]((1/2)M^1 - \cC M^1)$ in (\ref{Semicommutator}) is a pseudodifferential operator of order $-2$ on the surface $S$.
\hfill $\square$

One deduces from the proof of Lemma \ref{LemmaSC} that $T_{M_2 M_1} = T_{M_2} T_{M_1}$ holds actually up to trace class operators, if $n = 1$. 
In this case the results of \cite{HeltHowe75} apply to evaluate the index of Fredholm Toeplitz operators.

Lemma \ref{LemmaSC} allows one to develop the Fredholm theory of Toeplitz operators with operator-valued multipliers of the form (\ref{Equ_OM}).

\bigskip
\textit{The research of the author was supported by the Deutscher Akademischer Austauschdienst and by the Russian Federation President grant of support of leading scientific schools NSH-9149.2016.1.}

\end{document}